\documentclass[12pt]{article}
\usepackage{amssymb}
\font\smallit=cmti10
\font\smalltt=cmtt10
\font\smallrm=cmr9

\usepackage[usenames]{color}
\usepackage[colorlinks=true,
linkcolor=webgreen,
filecolor=webbrown,
citecolor=webgreen,
bookmarks=false]{hyperref}

\definecolor{webgreen}{rgb}{0,.5,0}
\definecolor{webbrown}{rgb}{.6,0,0}

\usepackage{epsfig}
\usepackage{amssymb}
\usepackage{amsmath}

\newcommand{\crit}{\alpha}

\newenvironment{packed_enumerate}{
\setlength{\parsep}{0pt}
\setlength{\parskip}{0pt}
\begin{enumerate}
  \setlength{\itemsep}{1pt}
  \setlength{\parsep}{0pt}
  \setlength{\parskip}{0pt}
}{\end{enumerate}}

\begin{document} 

\begin{center}
\vspace*{-20pt} 
\centerline{\smalltt INTEGERS: \smallrm ELECTRONIC JOURNAL OF COMBINATORIAL NUMBER THEORY
\smalltt 9 (2009), \#G01} 
\vskip 20pt

\uppercase{\bf The shortest game of Chinese Checkers \\ and related problems}
\vskip 20pt
{\bf George I. Bell}\\
{\smallit Tech-X Corporation, 5621 Arapahoe Ave Suite A, Boulder, CO 80303, USA}\\
{\tt gibell@comcast.net}\\
\end{center}

\vskip 30pt
\centerline{\smallit Received: 3/4/08, Revised: 10/13/08, Accepted: 12/20/08, Published: 1/5/09}
\vskip 30pt

\centerline{\bf Abstract}

\noindent
In 1979, David Fabian found a complete game of
two-person Chinese Checkers in 30 moves (15 by each player)
[Martin Gardner, Penrose Tiles to Trapdoor Ciphers, MAA, 1997].
This solution requires that the two players cooperate to
generate a win as quickly as possible for one of them.
We show, using computational search techniques, that no
shorter game is possible.
We also consider a solitaire version of Chinese Checkers where
one player attempts to move her pieces across
the board in as few moves as possible.
In 1971, Octave Levenspiel
found a solution in 27~moves [Ibid.];
we demonstrate that no shorter solution exists.
To show optimality, we employ
a variant of A* search,
as well as bidirectional search.

\pagestyle{myheadings}
\markright{\smalltt INTEGERS: \smallrm ELECTRONIC JOURNAL OF COMBINATORIAL NUMBER THEORY \smalltt 9 (2009), \#G01\hfill}

\thispagestyle{empty} 
\baselineskip=15pt 
\vskip 30pt

\section*{\normalsize 1. Introduction}

The game of \textbf{Halma} was invented in the 1880's by
George H. Monks \cite{MonksPatent, Whitehill}.
This game is played on a rather large $16\times 16$ board
and is still popular in parts of Europe.
In 1892, a significant variation appeared in Germany played on a triangular grid,
originally called \textbf{Stern-Halma} \cite{Whitehill}.
When this game was marketed in the United States it was given
the more exotic-sounding name \textbf{Chinese Checkers},
although it did not originate in China and is not a variant of Checkers.
Chinese Checkers remains a popular children's game in the United States.

\begin{figure}[htbp]
\centering
\epsfig{file=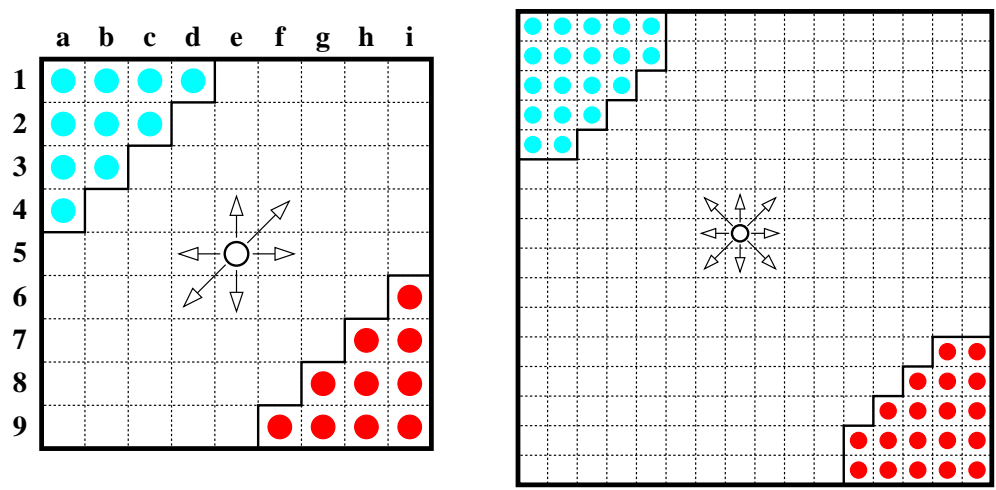}
\caption{A Chinese Checkers board (left) and Halma board (right)
with men in their starting positions.
The central man shows the directions of allowed movement.}
\label{fig1}
\end{figure}

Chinese Checkers is normally played using marbles on a 121-hole, star-shaped board.
However, the two-player version of the game can equivalently be played
on a square $9\times 9$ board, with move directions as given in
Figure~\ref{fig1}a\footnote{This is possible because
the two players move between two opposing ``points" of the star.
The 3 or 6-player versions of Chinese Checkers cannot be
played on a square $9\times 9$ board.}.
Although the board symmetry is harder to see,
it allows us to consider Halma and Chinese Checkers
as games played on the same board shape.

We will refer to a board location as a \textbf{cell},
with coordinates given in Figure~\ref{fig1}a.
We will also find it useful to use Cartesian Coordinates
to refer to a cell, with the origin at the center
of the board\footnote{The Halma board has no central hole,
we just use the cell identified in Figure~\ref{fig1}b as the origin.}.
Each player begins with a certain number of identical game pieces,
called \textbf{men}\footnote{To compensate for this male-oriented terminology,
all players will be female.}.
The set of men owned by one player will be referred to as her \textbf{army}.
The standard Chinese Checkers army has 10 men, and the standard Halma
army has 19 men.

A \textbf{move} in either game is
a step move or a jump move, as defined by:
\begin{packed_enumerate}
\item A \textbf{step}, where a man is moved one cell in the
direction indicated by the arrows in Figure~\ref{fig1}.
\item A \textbf{chain} of one or more jumps by the same man.
A \textbf{jump} is where a man hops over another man (of either army)
into an empty cell.
Jumps are allowed in any of the directions indicated by the arrows in Figure~\ref{fig1}.
The jumped piece is \textit{not removed} from the board---there
are no captures.
Jumps are never compulsory---a player may choose to stop a chain
whenever she pleases.
\end{packed_enumerate}

The area where a player's army begins from (outlined in Figure~\ref{fig1})
will be called that player's \textbf{base}.
Players alternate moves; the first player to fully occupy the
opposing base is the winner.
A complicating factor is how to deal with a player who refuses to
vacate their own base, preventing the other player from winning.
Various additional rules have been proposed to prevent this \cite{Whitehill},
but these will not concern us here.

In this square-board geometry the movement rules for Chinese Checkers and Halma
are identical except for the directions of possible movement.
Halma allows moves in all 8 directions from a cell, along rows, columns,
and both diagonals---we will refer to these movement rules
as \textbf{8-move} rules.
In Chinese Checkers, one parallel direction of diagonal movement is
not allowed, these will be referred to as \textbf{6-move} rules.
Although it may appear artificial to remove only one direction of diagonal jump,
this rule variation is completely natural on the triangular grid
on which Chinese Checkers is normally played.

\begin{figure}[htbp]
\centering
\epsfig{file=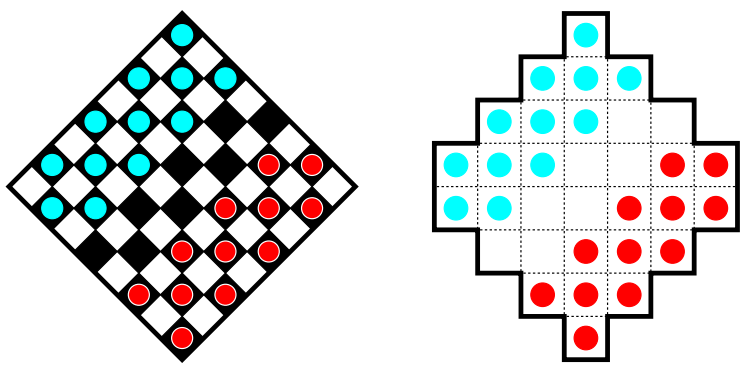}
\caption{(a) A standard game of checkers.
(b) The same game under 4-move rules.}
\label{fig2}
\end{figure}

A final variation is to allow steps and jumps only along
columns and rows;
this variation will be called \textbf{4-move} rules.
Figure~\ref{fig2}a shows a standard game of Checkers
(English Draughts) viewed with the board rotated $45^\circ$.
We can remove the white squares (they are not used),
and play Checkers on the 32-cell board in Figure~\ref{fig2}b
under 4-move rules.
Although Checkers has more complex jumping rules (captures),
this shows how the basic jumps and steps of the game can be viewed
as moves restricted to columns and rows---that is, 4-move rules.

Halma game play naturally divides into three distinct phases
(using George Monks' original terminology \cite{MonksPatent}):
\begin{packed_enumerate}
\item the \textbf{gambit}, where the armies advance toward one another but do not interact,
\item the \textbf{melee}, where the two armies interact and eventually pass through one another,
\item the \textbf{packing}, where the armies separate and attempt to fill the opposing bases as quickly as possible.
\end{packed_enumerate}

This paper does not address game play directly,
but considers two puzzles based on Chinese Checkers and Halma.
The first puzzle concerns the shortest possible game.
In 1979, David Fabian, working by hand, found a complete game
of Chinese Checkers in 30~moves \cite[p. 309]{GardnerPT}
(15~moves by each player).
This is remarkably short considering that each of the 20 men
(both armies) must make at least one move.
Such a solution requires that both sides cooperate so
that one of them wins as quickly as possible,
and thus has little to do with a competitive game.
We will show that no game can be shorter than 30~moves.

A second puzzle we consider is a solitaire version of the game,
where the goal is to advance an army across the board in as
few moves as possible (with no opponent's pieces in the way).
We'll refer to these puzzles as ``army transfer problems".
Such problems were a favorite of Martin Gardner---in three
of his Scientific American columns he
discusses army transfer problems on three different boards.
First with regard to a Checkers board (Figure~\ref{fig2})
under 4-move rules \cite{GardnerNMD}, then in the context of Halma
on a $9\times 9$ board under 8-move rules \cite{GardnerWL},
and finally in the context of Chinese Checkers
under 6-move rules \cite{GardnerPT}.
We will find the shortest possible solution to these three
army transfer problems, and consider various generalizations.

\vskip 30pt 
\section*{\normalsize 2. The shortest possible game}

The game is considered ended when one player fully occupies
her opponent's base, even if the other player does the same
on her next move.
Draws can't occur except in artificial situations where a player
refuses to leave her base, or the players move pieces back and
forth so that the game goes on forever.
On some boards one player can be stalemated with
no possible moves.
This is not possible in standard Chinese Checkers or Halma
(there are not enough men in one army to trap the opposing army in a corner),
but can occur among cooperating players for other variants.
We'll not consider such stalemated games as candidates for
the shortest possible game.

The length of a completed game is defined as the total number of moves
taken by both players.
We assume that the \textit{blue} player, starting from the upper left in
Figure~\ref{fig1}, always moves first.
The \textit{red} player moves second; if she wins first on her
15th move, it is a 30~move game.
If the blue player wins first on her 15th move,
it is a 29~move game.

\vskip 30pt 
\subsection*{\normalsize 2.1 Types of men}

\begin{figure}[htbp]
\centering
\epsfig{file=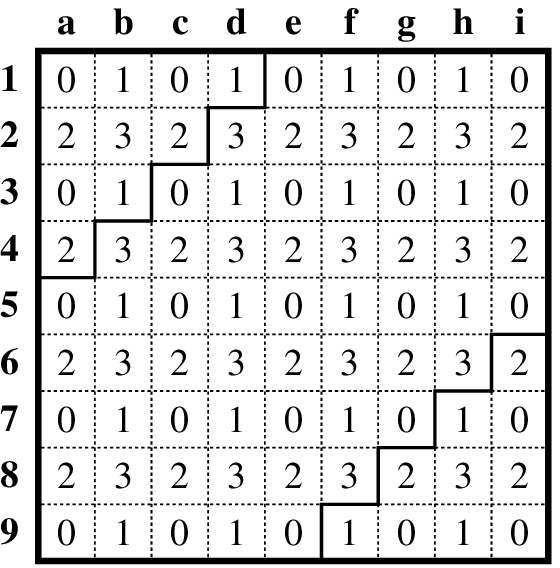}
\caption{Type labeling of a Chinese Checker board.}
\label{fig3}
\end{figure}

Clearly jumps are more effective than steps to move an army
quickly across the board.
However, a certain number of steps are generally needed.
A skilled player selects a careful mixture of jumps and steps
to advance her army across the board.

There is a fundamental difference between steps and jumps.
In Figure~\ref{fig3}, we checker the board with a pattern
of four type labels $0$--$3$.
We label the men in an army by their type, and observe that
\textit{only step moves can change a man's type}.
On a Chinese Checker board,
both bases have the same number of men of each type.
Therefore, the winning player will have the same types
of men at the end of the game that she started with.
Interestingly, this is not true of Halma, where the
type 0 and type 3 men must effectively change places
during the game.
The counts of the number of types of men for each game
are given in Table~\ref{table1}.

\begin{table}[htbp]
\begin{center} 
\begin{tabular}{ | c | c | c | c | c | }
\hline
 & \multicolumn{2}{| c |}{Chinese Checkers} &  \multicolumn{2}{| c |}{Halma} \\		
Type & \# starting & \# finishing &  \# starting & \# finishing \\
\hline
\hline
0 & 3 & 3 & 6 & 3 \\
\hline
1 & 3 & 3 & 5 & 5 \\
\hline
2 & 3 & 3 & 5 & 5 \\
\hline
3 & 1 & 1 & 3 & 6 \\
\hline
\hline
total men & 10 & 10 & 19 & 19 \\
\hline
\end{tabular}
\caption{The number of men of each type for the games of Chinese Checkers and Halma
(from the blue player's perspective).} 
\label{table1}
\end{center} 
\end{table}

This suggests that it may be possible to play an entire game of Chinese Checkers
without making any step moves, and we will see that it is possible.
For Halma, however, the same argument shows that
\textit{at least three step moves are required} to win the game.
The difference is due to the fact that the board side in Halma is
even (20), while for Chinese Checkers it is odd (9).
In general, we'll refer to these as \textbf{even} or \textbf{odd} boards,
and this type disparity will be seen on any even board.
On an odd board, the number of starting and finishing
types will be the same, as long as the starting army is symmetric
about the diagonal line $x=-y$.

In Chinese Checkers, to advance a type 0 or 3 man (move it closer to the
opponent's base) requires a jump over a type 1 or 2 man,
or a step, which converts it to a type 1 or 2 man.
Similarly, advancing a type 1 or 2 man requires a jump over or conversion
to a type 0 or 3 man.
In Halma, a man of a certain type can advance by jumping or converting into
any of the other three types
(because it is now possible to move diagonally between the bases).

\vskip 30pt 
\subsection*{\normalsize 2.2 Ladders}

\begin{figure}[htbp]
\centering
\epsfig{file=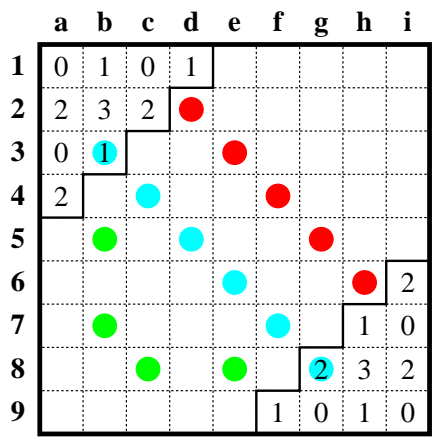}
\caption{Three ladders on a Chinese Checkers board.}
\label{fig4}
\end{figure}

A \textbf{ladder} is any
configuration of men that
allows for quick transfer of men between the bases by means
of long chain jump moves.
Figure~\ref{fig4} shows several ladders on a Chinese Checkers board.
These ladders are shown using the same color men to visually separate them,
but in general they can be composed of either army.

The red (top) ladder in Figure~\ref{fig4} is composed of men of types 0 and 3,
and can transport men of the opposite types 1 and 2.
We will refer to such a ladder as a ``type 0\&3 ladder".
The blue (middle) ladder has these roles reversed, it is a type 1\&2 ladder.
The green (bottom) ladder is also a type 1\&2 ladder,
but can only transport type 3 men.
Because it can transport only one type of man, the green ladder is less useful.
Any ladder can transport at most two of the four types of men,
for the fastest transport of an army (with men of all types),
\textit{at least two ladders of different types will be needed}.

The game would be relatively easy if the ladders were in place at the
beginning.
In reality, ladders are not only used, they must be built and disassembled.
Plus, in a competitive game ladders will seldom be separate
or complete as in Figure~\ref{fig4},
and one player can strategically ``block" a ladder by stopping
a man in the middle.
But our focus now is on the shortest possible game,
where the players must cooperate to build the most efficient ladders,
and they need not be broken down at all if they are composed
of the losing player's men.

\vskip 30pt 
\subsection*{\normalsize 2.3 A lower bound on game length}

Given two cells with Cartesian coordinates
$\mathbf{a}=(a_x,a_y)$ and $\mathbf{b}=(b_x,b_y)$, what is the
minimum number of step moves needed to move a man between the two?
To answer this question, we define the norm of a cell
with Cartesian coordinates $(x,y)$ as
\begin{equation}
||(x,y)|| =\left\{ \begin{array}{lll}
||(x,y)||_{\infty} &= \max(|x|,|y|) &  \mbox{Halma or 8-move,} \\
||(x,y)||_{1} &= |x| + |y| & \mbox{4-move,} \\
||(x,y)||_{\vartriangle} &= \frac{1}{2}(|x| + |y| + |x-y|) & \mbox{Chinese Checker or 6-move.}
\end{array} \right.
\label{eq:normdef}
\end{equation}
The norm $||(x,y)||_{\vartriangle}$ is a combination of the first two norms,
as shown by the alternate formula
\begin{align}
||(x,y)||_{\vartriangle} &=
\left\{ \begin{array}{ll}
\max(|x|,|y|) &  \mbox{if } \mbox{sgn}(x)=\mbox{sgn}(y), \\
|x| + |y| & \mbox{otherwise.} \\
\end{array} \right. \label{eq:altnormdef}
\end{align}
The formula (\ref{eq:altnormdef}) was given in 1976 in the context
of image processing to calculate distances on a hexagonal grid \cite{LuczakRos76}.

The distance between two board locations
$\mathbf{a}=(a_x,a_y)$ and $\mathbf{b}=(b_x,b_y)$ is then defined as
\begin{equation}
d(\mathbf{a},\mathbf{b}) = ||\mathbf{a}-\mathbf{b}|| = ||(a_x-b_x, a_y-b_y)|| ,
\label{eq:ptdist}
\end{equation}
with norm appropriate for the game as defined in (\ref{eq:normdef}).
$d(\mathbf{a},\mathbf{b})$ is the minimum number of steps
needed to move a man from $\mathbf{a}$ to $\mathbf{b}$.

Given two armies $B$ and $R$, we define the
distance between them as:
\begin{equation}
d(B, R) = \min\{ d(\mathbf{b}, \mathbf{r}), \forall \mathbf{b} \in B \mbox{ and } \mathbf{r} \in R \} .
\label{eq:armydist}
\end{equation}
If we let $B$ and $R$ be the initial positions of two armies
at the beginning of a game, for Chinese Checkers we
can compute that $d(B,R)=10$.
Interestingly, for Halma we also have $d(B,R)=10$.
Thus, despite the large difference in board size,
the initial distance between the two armies is the same for
Chinese Checkers and Halma.

{\bf Theorem 1} If $B$ and $R$ are the initial positions of
two equal-sized armies ($s=|B|=|R|$),
no game can be shorter than $h$ moves, where
\begin{equation}
h = \max\{0, d(B, R)-2\} + 2s - 1 .
\label{eq:gamelb}
\end{equation}

{\it Proof}: After one player moves, $d(B,R)$ can decrease by at most $1$.
As soon as $d(B,R)\le 2$, the next move can cross between the two armies.
The best that can happen is that $d(B,R)$ decreases to $2$,
the next player then wins by placing a man from her army into her
opponent's base on \textit{each} of her subsequent $s$ moves.
The number of moves is then given by (\ref{eq:gamelb}).

For a standard game of Chinese Checkers, the bound given by Theorem 1
is 27~moves, and for Halma 45~moves.

\vskip 30pt
\subsection*{\normalsize 2.4 A 30~move game of Chinese Checkers}

David Fabian found a 30~move Chinese Checkers game
in 1979 and sent it to Martin Gardner \cite[p. 309]{GardnerPT}.
I contacted David Fabian about how he found his solution,
he said he used ``logic and patience" to find it by hand.
We can come up with a set of properties that a short solution
is likely to have.
We can't prove that the shortest solution must have these properties,
but they can guide us on our search for them,
both by hand and using a computer.

From Section~2.1,
if the winning player in Chinese Checkers makes one step,
she will be obliged at some point to make another.
Therefore, it seems reasonable that in the shortest possible game,
the winning player only jumps.
The losing player has no constraints on the types of finishing men,
and is likely to make some steps.

If the game is to have length $L$, then there is a critical move,
$\crit=L-2(|A|-1)$, after which the winning player \textit{must} have
at least one man in her opponent's base.
In Section~2.2 we learned that for a quick game two ladders
must be built.
The first ladder can be built by both players in the moves before $\crit$,
while the second ladder must be finished by the losing player
in the moves after $\crit$.

\begin{figure}[htb]
\centering
\epsfig{file=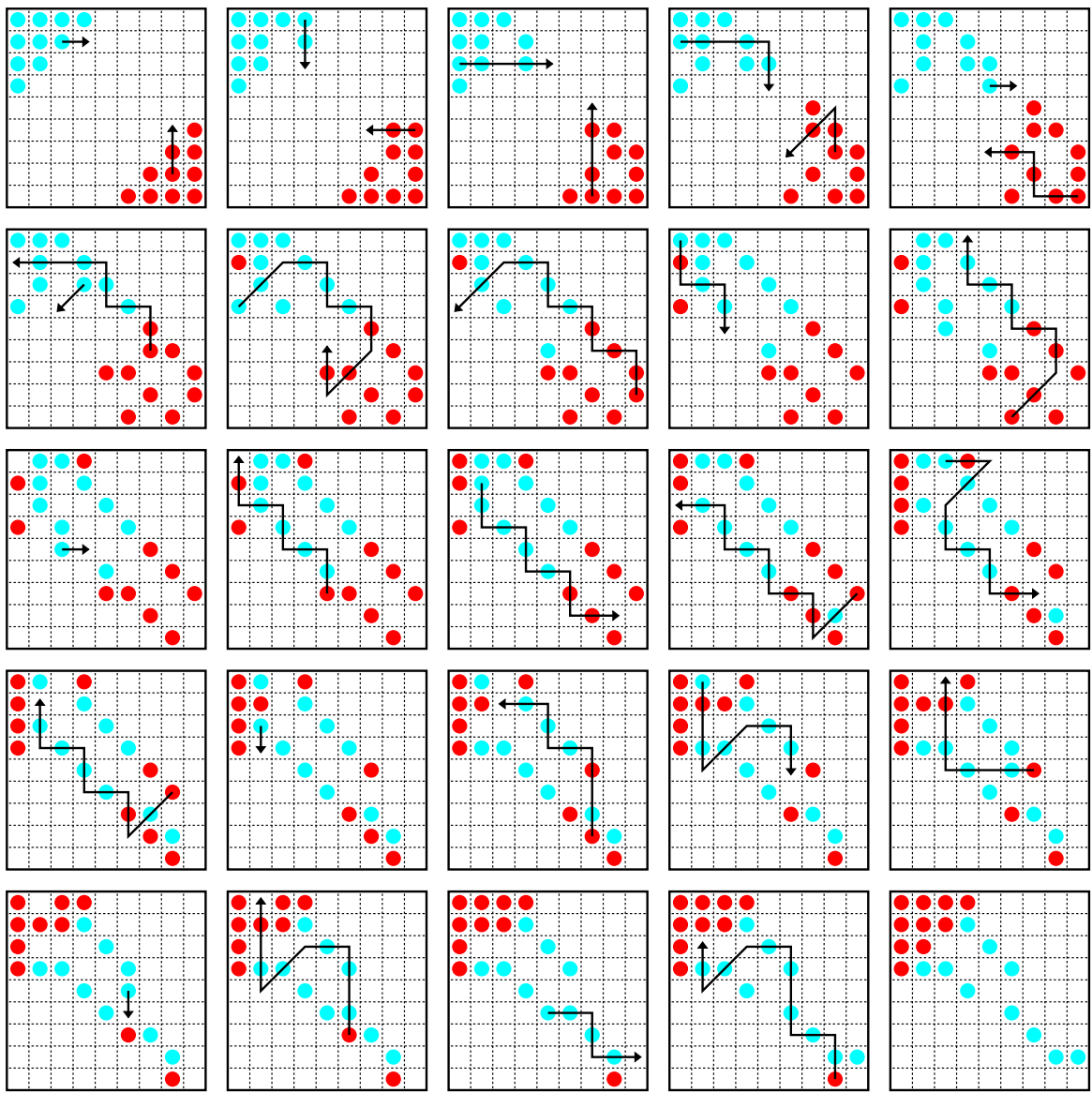}
\caption{David Fabian's 30 move game of Chinese Checkers.}
\label{fig5}
\end{figure}

In summary, some of the properties that a short solution usually has are:
\begin{packed_enumerate}
\item The winning player only jumps (on odd boards).
\item The first ladder is built by both players
during the first $\crit$ moves.
\item The second ladder is completed by the losing player on moves after $\crit$.
\item After each move $n\ge \crit$ by the winning player,
there must be at least $(n-\crit)/2+1$ men in the opposing base.
\item In the middle of the game, each diagonal line $x-y=k$ between the bases
must be occupied by at least one man,
where $-4\le k \le 4$ (this forces placement of two ladders).
\end{packed_enumerate}

Finally, Figure~\ref{fig5} shows David Fabian's 30~move solution,
which has all these properties ($\alpha=12$ in this case).
We have modified the losing player's last move slightly,
as there are many possibilities.
In fact this freedom might suggest the existence of a $29$
move solution, but we shall see that this is not the case.

\vskip 30pt
\subsection*{\normalsize 2.5 Search algorithms}

One way to show that no 29~move solution exists
is to do a search of the game tree.
We perform a breadth-first search so that we can easily eliminate
duplicate board positions.
The search proceeds by \textbf{levels},
where the level set $L_i$ consists of board positions
that can be reached after $i$ moves.
The level set $L_0$ contains only the initial board position,
while $L_1$ contains one element
for each possible first move by the blue player.
There are 14 possible first moves in Chinese Checkers,
but 7 of these are equivalent by symmetry,
so $|L_1|=7$, $|L_2|=7\cdot 14=98$,
and $|L_3|=1253$.
The total size of the state space for the game is
$\binom{81}{10,10}\approx 8.67\times 10^{23}$,
and the level sets $L_i$ grow much too rapidly to calculate
$L_{29}$.

Fortunately a very good lower bound on the number of
moves remaining exists, namely the bound given in Theorem~1.
Although the bound in Theorem~1 assumes the board is in the
initial position, it can be easily modified to apply to any
board position.
If the distance between the two armies $d(R,B)>2$,
we can use (\ref{eq:gamelb}) unmodified.
When this distance becomes two or less,
we modify $s$ in (\ref{eq:gamelb}) to be the number of men in
the winning army that are not in the opposing base.
Given any board position $P$, we can use (\ref{eq:gamelb})
as a lower bound on the number of moves remaining, $h(P)$.

To apply this in our search scheme we use a version of A* search
called ``breadth-first iterative deepening A*" \cite{ZhouAI}.
Suppose we are searching for a solution with length at least $m$.
At move $i$ if the board position is $P$,
no solution from this board position can be shorter
than $i+h(P)$.
Thus, we can terminate the search from this node
if $i+h(P)>m$.

This gives a search tree that expands rapidly until a critical move $\crit$,
the first move where the winning player must place at least one man in
the opposing base.
At this level the search tree contracts significantly.
But the next move by the losing player is
unconstrained so the search expands,
only to contract on the next move by the winning player
who again must place another man in the opposing base.

For Chinese Checkers, we know that a solution of length $30$ exists.
Therefore, we apply the search algorithm to look for all solutions
of length $29$ and $28$.
These searches come up empty, so the shortest solution has $30$ moves.
The solution in Figure~\ref{fig5} is not unique;
the search technique finds several hundred different
30~move games.

The same algorithm can find the shortest game for many other
starting configurations.
The search strategy in these cases proceeds follows:
first, by incorporating the heuristic rules of Section~2.4,
we find a solution of length $N$ which we believe is
the shortest possible.
Then, we remove these heuristic rules and use only the bound on
solution length (\ref{eq:gamelb}) to show that no
solution exists of length $m=N-1$ and $m=N-2$.

\begin{table}[htbp]
\begin{center} 
\begin{tabular}{ | l | c | c | c | c | c | c | }
\hline
 & & Move & Army & & Lower & Shortest \\
Game & Board & Rule & Size & $d(B,R)$ & Bound (\ref{eq:gamelb}) & Game \\
\hline
\hline
Chinese Checkers & $9\times 9$ & 6-move & 10 & 10 & $27$ & $30$ \\
\hline
Chinese Checkers (15 man) & $9\times 9$ & 6-move &  15 & 8 & $35$ & $36$ \\
\hline
Halma & $16\times 16$ & 8-move & 19 & 10 & $45$ & Unknown \\
\hline
Grasshopper & $8\times 8$ & 8-move & 10 & 4 & $21$ & $24$ \\
\hline
Chinese Checkers (8-move) & $9\times 9$ & 8-move & 10 & 5 & $22$ & $24$ \\
\hline
\end{tabular}
\caption{Summary of shortest game lengths on various boards.} 
\label{table2}
\end{center} 
\end{table}

Table~\ref{table2} shows the results of these runs.
For a Chinese Checkers game with 15~men per side,
David Fabian found a 38~move game \cite[p. 309]{GardnerPT},
but our algorithm finds that a 36~move game is the shortest possible.
Grasshopper \cite[p. 117]{GardnerWL} is Halma played
on an $8\times 8$ (even) board, with the initial configuration
of 10 men the same shape as in Chinese Checkers.
By the type analysis in Section~2.1, the winning player must make at
least 2 step moves, and a game in 22 moves or less is impossible.
In order to win in exactly 24~moves, the winning player must
make 2 step moves and 10 jump moves,
and each of the jump moves must finish either inside their
opponents base or one cell short of it.
For the same game on a $9\times 9$ board (odd),
the shortest game also has 24~moves,
but the winning player only jumps.

\vskip 30pt
\section*{\normalsize 3. Army transfer problems}

We now consider the problem of moving one army quickly between bases,
without any opposing men.
On a Chinese Checkers board,
this problem was discussed in a 1976 Martin Gardner column \cite{GardnerPT},
although the earliest reference is a 1959 Canadian periodical
on magic \cite{Ibidem}.
Octave Levenspiel worked by hand to find short transfers;
in 1971 he found a 27~move solution \cite[p. 73]{GardnerPT}.

\begin{figure}[htb]
\centering
\epsfig{file=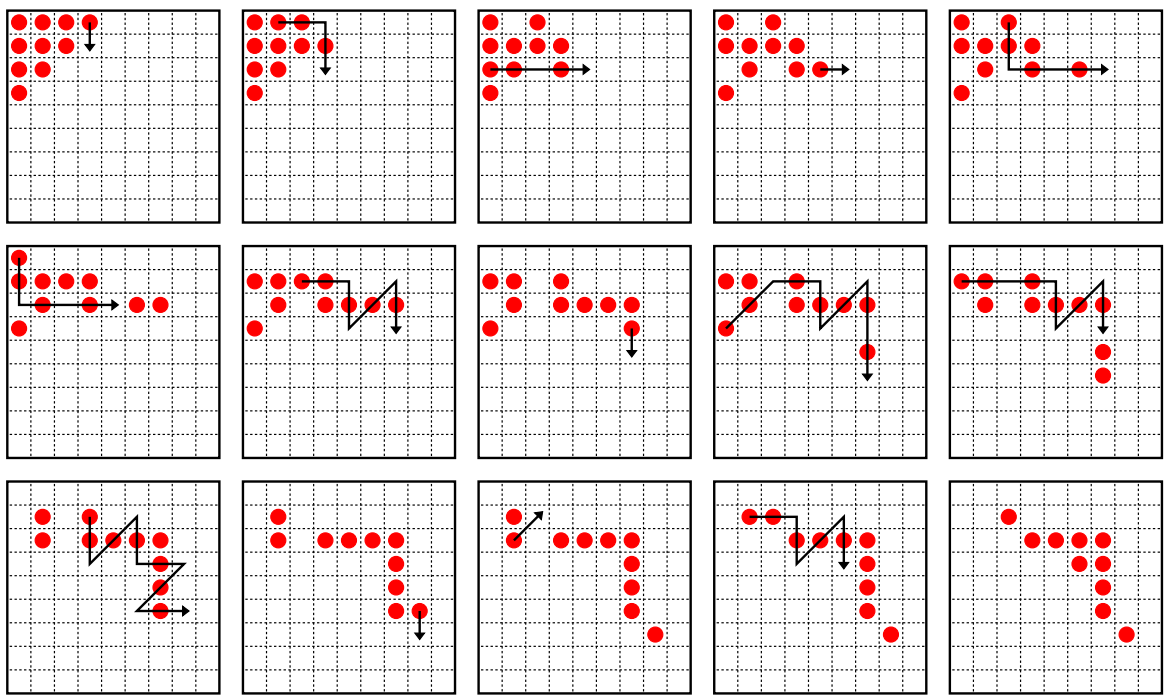}
\caption{Octave Levenspiel's 27 move solution.
Only the first half of the solution is shown, because the solution is a palindrome
(only half of the 14th move is shown).}
\label{fig6}
\end{figure}

Figure~\ref{fig6} shows Levenspiel's 27~move solution.
An interesting feature of this solution is that it is \textbf{palindromic},
meaning that
moves 15--27 are mirror images of the moves 1--13,
taken in the reverse order
(where the mirror symmetry is about the line $x=y$),
and move 14 is itself mirror symmetric.
If the solution is interrupted in the middle of the
14th move (the last diagram of Figure~\ref{fig6}),
the board position is mirror symmetric about the line $x=y$.
It is \textit{not true} that every 27~move solution to this problem is palindromic,
but palindromic solutions often seem to exist,
and are the most elegant.

In 1973, Harry O. Davis sent Martin Gardner a proof that 27~moves was the
shortest possible solution to this problem.
Although this proof is mentioned in Gardner's book \cite[p. 68]{GardnerPT},
it was never published.
This proof has been preserved in Martin Gardner's files \cite{Davis}.
Davis begins his one-page proof with an argument that
the shortest solution to the problem must contain at least 10 step moves
(note that Levenspiel's solution in Figure~\ref{fig6}
contains 10 steps and 17 jump moves).
However, we have found 27~move solutions to this problem
with only 8 step moves (see Appendix A),
so Davis' claim is false.
Although the proof appears flawed,
the theorem is nonetheless true---there is no solution shorter than 27~moves.
This will be demonstrated numerically below.

\vskip 30pt
\subsection*{\normalsize 3.1 Centroids and bounds}

The (diagonal) \textbf{centroid} of a man $\mathbf{a}$
with Cartesian coordinates $(a_x,a_y)$ is defined as
\begin{equation}
c(\mathbf{a})=a_x-a_y.
\label{eq:mcentdef}
\end{equation}
The centroid $c(\mathbf{a})$ is a measure of how far this man has
progressed in his journey between the bases.
In the blue player's army,
the centroid of any man begins at $-5$ or less,
and ends at $+5$ or greater.

The \textbf{centroid} of an army $A$ is the integer-valued function defined by
\begin{equation}
c(A) = \sum_{\mathbf{a}\in A} c(\mathbf{a}) = \sum_{\mathbf{a}\in A} a_x-a_y.
\label{eq:acentdef}
\end{equation}
This is related to the classical centroid, or center of gravity,
in the following sense:
the army will balance across the diagonal line
$x-y=c(A)/|A|$,
where $|A|$ is the size of army $A$.

This centroid gives a natural measure of an army's progress,
and in competitive games which player is currently in the lead.
In a game of Chinese Checkers, for example,
where the origin is the center of the board,
for the blue player
$c(A)$ begins at $-60$ and ends at $+60$.
Note that the starting army balances on the line $x-y=-6$.
When $c(A)=0$, the army is exactly half-way to their goal,
and balances across the diagonal line $x=y$.

{\bf Theorem 2} Consider an army $A$, and let $\mathbf{a}_{min}\in A$ be
a man with minimum centroid,
i.e. $c(\mathbf{a}_{min})\le c(\mathbf{a}), \forall \mathbf{a}\in A$.
Similarly, let $\mathbf{a}_{max}\in A$ a man with maximum centroid.
Then in one move, the centroid can increase at most
$\delta$, or decrease at most $-\delta$, where
\begin{equation}
\delta  = c(\mathbf{a}_{max})-c(\mathbf{a}_{min})+\ell,
\label{eq:ubound}
\end{equation}
where $\ell=1$ for 4-move and 6-move play,
and $\ell=2$ for 8-move play.

{\it Proof}: Only one man can move,
and the greatest centroid increase is achieved by taking a man $\mathbf{a}_{min}$
with minimum centroid and increasing his centroid as much as possible.
The best $\mathbf{a}_{min}$ can do is finish with a jump over some man with
centroid $c(\mathbf{a}_{max})$.
Under 8-move play, if the last jump is diagonal, his final centroid
can be be at most $c(\mathbf{a}_{max})+2$, otherwise it can be at most $c(\mathbf{a}_{max})+1$.

Theorem~2 gives a relatively crude upper bound on centroid increase,
which it is tempting to refine further.
Under 4 and 6-move rules,
if $c(\mathbf{a}_{max})-c(\mathbf{a}_{min})$ is even,
in order to reach $c(\mathbf{a}_{max})+1$, the man $\mathbf{a}_{min}$
must finish with a rightward or downward jump over a man $\mathbf{a}_{max}$,
and this is not possible because this man has the wrong type.
In this case the upper bound (\ref{eq:ubound}) can be reduced by one.
More significantly, suppose the army can be partitioned into two pieces
$A_1$ and $A_2$ with $\mathbf{a}_{min}\in A_1$, $\mathbf{a}_{max}\in A_2$
and $d(A_1,A_2)>2$.
Then the man $\mathbf{a}_{min}$ cannot reach the other half of the army,
and the bound (\ref{eq:ubound}) can be reduced considerably.
All such refinements result in a more complex
formula for centroid increase.

We will also want to use Theorem~2 iteratively
to obtain an upper bound on the
centroid increase of an army after $n$ moves.
In this situation, the complexity of an improved bound increases greatly
because we will have to consider the interaction of multiple moves.
For example, in the case where the army can be partitioned into
$A_1$ and $A_2$ with $d(A_1,A_2)=3$,
the first move could connect the two pieces,
and then the second move could go between them.
The crude bound in Theorem~2 is much easier to implement iteratively,
because it is valid no matter how the moves interact.

We can also use the ideas in Theorem~2 to get a lower bound on the number
of moves to accomplish the transfer.
We assume that the initial army is $A$, and the location of the target base is $B$
(a fixed set of cells).
As the army $A$ advances, the distance $d(A,B)$ can decrease by at most one per move.
Thus, the smallest number of moves to place a man in the target base $B$ is $d(A,B)$.
After this, the remaining $|A|-1$ men move into the base $B$ at best one per move,
so any solution to the transfer problem has length at least $d(A,B)+|A|-1$.
For the standard Chinese Checkers army, this gives a lower bound of 19~moves,
not a very tight bound considering the minimum is 27 (as we shall see).
Although these crude estimates don't get very close to the true optimum,
we will see that Theorem~2 is useful during a numerical search.

\begin{figure}[htb]
\centering
\epsfig{file=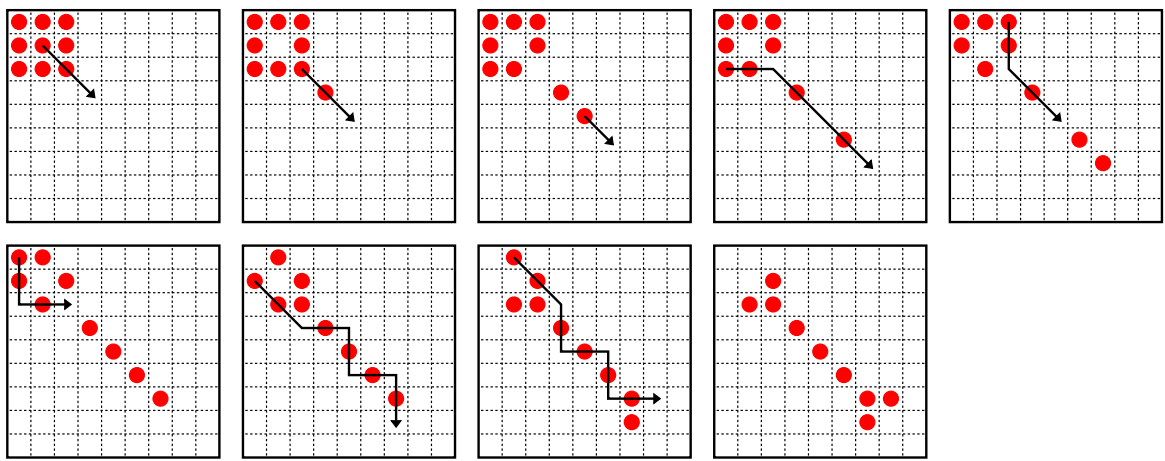}
\caption{An elegant 16~move solution to a Halma transfer problem
(H. Ajisawa and T. Maruyama \cite[p. 118]{GardnerWL}).
Again, only the first half of the solution is shown.}
\label{fig7}
\end{figure}

Figure~\ref{fig7} shows a 16~move solution to a Halma transfer problem
\cite[p. 118]{GardnerWL}.
The starting and ending armies in this case are $3\times 3$ square arrays of men.
Diagonal steps and jumps are allowed,
in our terminology we have 8-move rules.
The solution has even length and is also palindromic,
the board position after 8~moves is mirror symmetric about the line $x=y$.
Although our simple lower bound indicates that any solution must be at least 12~moves,
we'll be able to show using the search algorithm in the next section
that the 16~move solution in Figure~\ref{fig7} is the shortest possible.

\vskip 30pt
\subsection*{\normalsize 3.2 Search algorithms}

As before our basic computational tool is breadth-first search, with
the level sets $L_i$ defined as before, except that only one player moves.
Since the centroid must increase by $120$ (from $-60$ to $+60$),
clearly the most productive moves are those that increase the centroid.
To cut down on the size of the level sets, it would seem we should
\textit{only} consider moves that increase the centroid.
If we do so, however, there is no guarantee that
we will find the shortest possible solution.
Surprisingly, solutions of shortest possible length may contain
moves that \textit{decrease the centroid}.
Examples of this counter-intuitive phenomenon will be found below.

There are several symmetries in this problem.
The first is that the starting army is symmetric about the
diagonal line $x=-y$.
This effectively halves the number of possible starting moves,
and decreases the search space by a factor of two.
The starting and finishing positions are also symmetric about
the diagonal line $x=y$.
Thus, the set of possible board positions one move before the finish
is the same as the set of possible board positions one move
from the start, reflected across the line $x=y$.
This suggests that a good search technique is
\textbf{bidirectional search} \cite{bidirectional}.
To search for a solution of length $2N$,
we do a breadth-first search to $N$ moves, and intersect
the level set $L_N$ with the set obtained by reflecting
each element of $L_N$ across the diagonal line $x=y$.

Suppose we have a solution of length $2N+1$, and want to prove that
no solution of length $2N$ exists.
Suppose we also know that the maximum possible centroid
after \textit{any} sequence of $N$ moves is $C^{max}_N$
(this is the maximum centroid of any board position in the level set $L_N$).
For a solution of length $2N$ to exist, it must be that $C^{max}_N\ge 0$.
The key observation is that \textit{any solution} of
length $2N$ must have centroid \textit{at least} $-C^{max}_N$ at level $N$.

At any level $i\le N$, we apply Theorem~2 iteratively to get an
upper bound on the centroid after $N-i$ additional moves.
If this upper bound is less than $-C^{max}_N$,
then this board position cannot lead to a solution and
we need not search further from this board position.
This search combines aspects of bidirectional and A* search techniques,
and will be called a minimum centroid constraint (MCC) search,
because it eliminates boards with centroid too small to lead to a solution.

If you try to implement this MCC search,
you will notice a problem.
We assumed that we knew $C^{max}_N$ at the start of the algorithm,
yet this number is not determined until the algorithm is finished!
In reality we must estimate $C^{max}_N$,
then run the search algorithm using this estimate.
If the search algorithm produces a board at level $N$ with centroid greater
than $C^{max}_N$, we must run it all over again with the correct $C^{max}_N$.
It is only after the search finishes with a self-consistent value of $C^{max}_N$,
and the bidirectional search comes up empty, that we are assured that
no solution of length $2N$ exists.

One way to estimate $C^{max}_N$ is to truncate the search at each level,
keeping only the top $M$ boards with the largest centroid.
We have found that an $M$ of a few million gives a good estimate of $C^{max}_N$.
Short solutions can also be found quickly using this truncation technique.
Bill Butler solved the Chinese Checker transfer problem (Figure~\ref{fig6}) this way,
as documented on his web site \cite{DurangoBill}.
He found five different 27~move solutions, but because he truncated the number of
boards at each level, his search is not exhaustive.

For the standard Chinese Checker transfer problem,
we find that $C^{max}_{13}=5$.
We should note that although the MCC algorithm makes the
search problem solvable, it is still not easy.
The largest level set $L_{10}$ contains $1.3\times 10^7$ boards,
and the complete (unsuccessful) search for a 26~move solution
takes over two hours of CPU time\footnote{On a 1 GHz PC with 512MB of RAM.}.

\begin{table}[htbp]
\begin{center} 
\begin{tabular}{ | c | l | c | c | c | c | c | }
\hline
 & & & & \multicolumn{3}{| c |}{Shortest solution under}  \\
Problem & Configuration & \# Men & Board & 4-move & 6-move & 8-move \\
\hline
\hline
\# 1 & Triangle & 6 & $9\times 9$ & 25 & 23 & 16 \\
\hline
\# 2 & Triangle & 10 & $9\times 9$ & 30 & 27 (Fig~\ref{fig6}) & 20 \\
\hline
\# 3 & Triangle (jumps only) & 10 & $9\times 9$ & 46 & 35 (Fig~\ref{fig8}) & 21 \\
\hline
\# 4 & Triangle & 15 & $9\times 9$ & $36\dagger$ & $31\dagger$ (Fig~\ref{fig9}) & $26\dagger$ \\
\hline
\# 5 & Square & 4 & $9\times 9$ & $15$ & 15 & 12 \\
\hline
\# 6 & Square & 9 & $9\times 9$ & 25 & 25 & 16 (Fig~\ref{fig7}) \\
\hline
\# 7 & Checkers start & 12 & Fig~\ref{fig2}b & 20 & 16 & 16 \\
\hline
\hline
\multicolumn{7}{| l |}{$\dagger$ - MCC search did not terminate, a shorter solution may be possible.} \\
\hline
\end{tabular}
\caption{Summary of shortest solution lengths for transfer problems, all are the shortest possible except as noted.}
\label{table3}
\end{center} 
\end{table}

Table~\ref{table3} gives the results of such a search over a wide variety of problems,
with solutions given in Appendix~A.
Except as noted, the MCC search is run to completion in each case,
so the number of moves given is the smallest possible.
In particular, Octave Levenspiel's 27~move solution in Figure~\ref{fig6} is the shortest possible,
and 16~moves is the shortest solution to the Halma transfer problem in Figure~\ref{fig7}.
Problem \#7 was suggested by Gardner \cite{GardnerNMD},
and refers to moving a 12-man army across a checkers board.
This is equivalent to moving one of the 12-man armies in
Figure~\ref{fig2}b to the opposite side under 4-move rules.
The shortest solution to this problem has 20~moves.

\begin{figure}[htb]
\centering
\epsfig{file=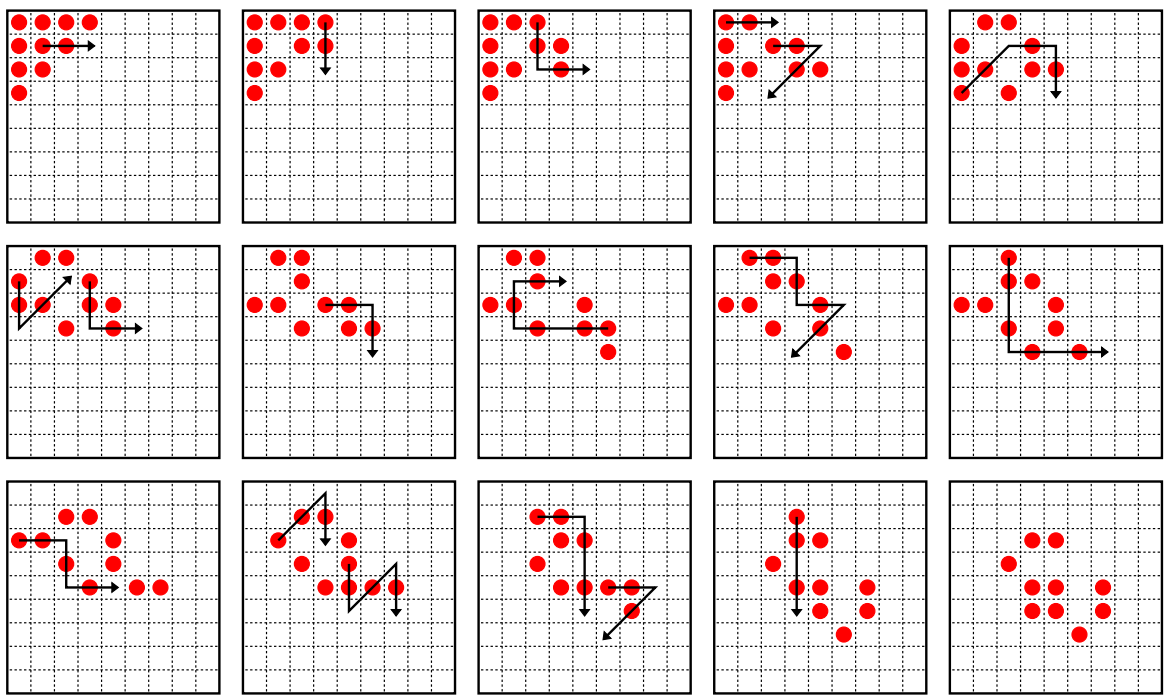}
\caption{The first half of a 35~move solution with jumps only
(the shortest possible).
Note that the 10th move (middle diagram, second row) goes backwards!}
\label{fig8}
\end{figure}

If we allow only jump moves, how quickly can a 10-man triangular army be transferred?
Table~\ref{table3} shows that the answer (under 6-move rules) is 35~moves,
with a sample solution shown in Figure~\ref{fig8}.
One interesting feature of this solution is a
\textit{backward} jump on the 10th move which
reduces the army's centroid $c(A)$.
This backward jump move can be explained by the fact that the man making it is of type 3.
There is only one type 3 man in this army, and in order to utilize him he must
be moved around a lot, even jumping backwards.
This phenomenon is also observed in the 4-move version of the problem.

Appendix~A includes solutions to problems in Table~\ref{table3} with
backward step moves.
These exceptions indicate that we cannot eliminate backward jump and step moves when searching
for the shortest possible solution.

One interesting problem is the 15-man Chinese Checkers army
(problem \#4 in Table~\ref{table3}).
In 1974, Min-Wen Du of Taiwan sent Martin Gardner
a letter with a solution to this problem
in 35 moves\footnote{Note that \cite[p. 68]{GardnerPT}
describes this problem incorrectly.
The letter \cite{Du} removes any ambiguity.} \cite{Du}.
My search algorithm has found a solution
to this problem in 31~moves (see Figure~\ref{fig9}),
but an MCC search doesn't finish,
so it is not known if this is the shortest solution possible.
This 31~move solution is interesting because
it uses two ladders rather than one.

\vskip 30pt
\subsection*{\normalsize 3.3 Symmetry and palindromes}

In many cases the shortest solution can be chosen to be a palindrome,
and these solutions are arguably the most elegant.
We can add further constraints to search specifically
for palindromic solutions.

First, we define the level of symmetry of an army $A$.
Given a man $\mathbf{a}$ with Cartesian coordinates $(a_x,a_y)$,
the coordinates of the cell reflected across the line $x=y$ are $(a_y,a_x)$.
We define the function
\begin{equation}
\mbox{sym}(\mathbf{a}) =\left\{ \begin{array}{lll}
+1 &  \mbox{if } (a_y,a_x) \mbox{ is occupied}, \\
0  & \mbox{otherwise}.
\end{array} \right.
\label{eq:symdef}
\end{equation}
We then define the \textbf{army symmetry} of an army $A$ as
\begin{equation}
\mbox{sym}(A)=\sum_{\mathbf{a}\in A} \mbox{sym}(\mathbf{a}),
\label{eq:armysymdef}
\end{equation}
so $\mbox{sym}(A)$ varies from zero to $|A|$ depending on the symmetry of the army across the line $x=y$,
and $\mbox{sym}(A)=|A|$ if and only if the army is mirror symmetric across this line.

Starting with an army $A_0$, denote by $A_i$ the board position after $i$ moves.
For a palindromic solution of odd length $2N+1$ to exist, it must be that
$\mbox{sym}(A_N)=\mbox{sym}(A_{N+1})=|A|-1$.
For a palindromic solution of even length $2N$ to exist,
we must have $\mbox{sym}(A_N)=|A|$.
These two situations can be seen in the second to the last diagram of
Figures~\ref{fig6} and \ref{fig8},
and the final diagram of Figure~\ref{fig7}.

The symmetry of the starting army, $\mbox{sym}(A_0)=0$, but by
the middle of a palindromic solution of length $2N$ or $2N+1$,
$\mbox{sym}(A_N)=T$,
where $T=|A|$ for an even solution length and $T=|A|-1$ for an odd solution length.
But the function $\mbox{sym}()$ cannot increase from $0$ to $T$ in one move.
Since a single move only affects one man,
$\mbox{sym}()$ can increase or decrease by at most two
(the only terms that can increase in the sum (\ref{eq:armysymdef})
are the endpoint of the move and its mirror image).
So for any $i\le N$, the board position $A_i$
in a palindromic solution must satisfy
\begin{equation}
\mbox{sym}(A_i)\ge T-2(N-i).
\label{eq:symconst}
\end{equation}

\begin{figure}[htb]
\centering
\epsfig{file=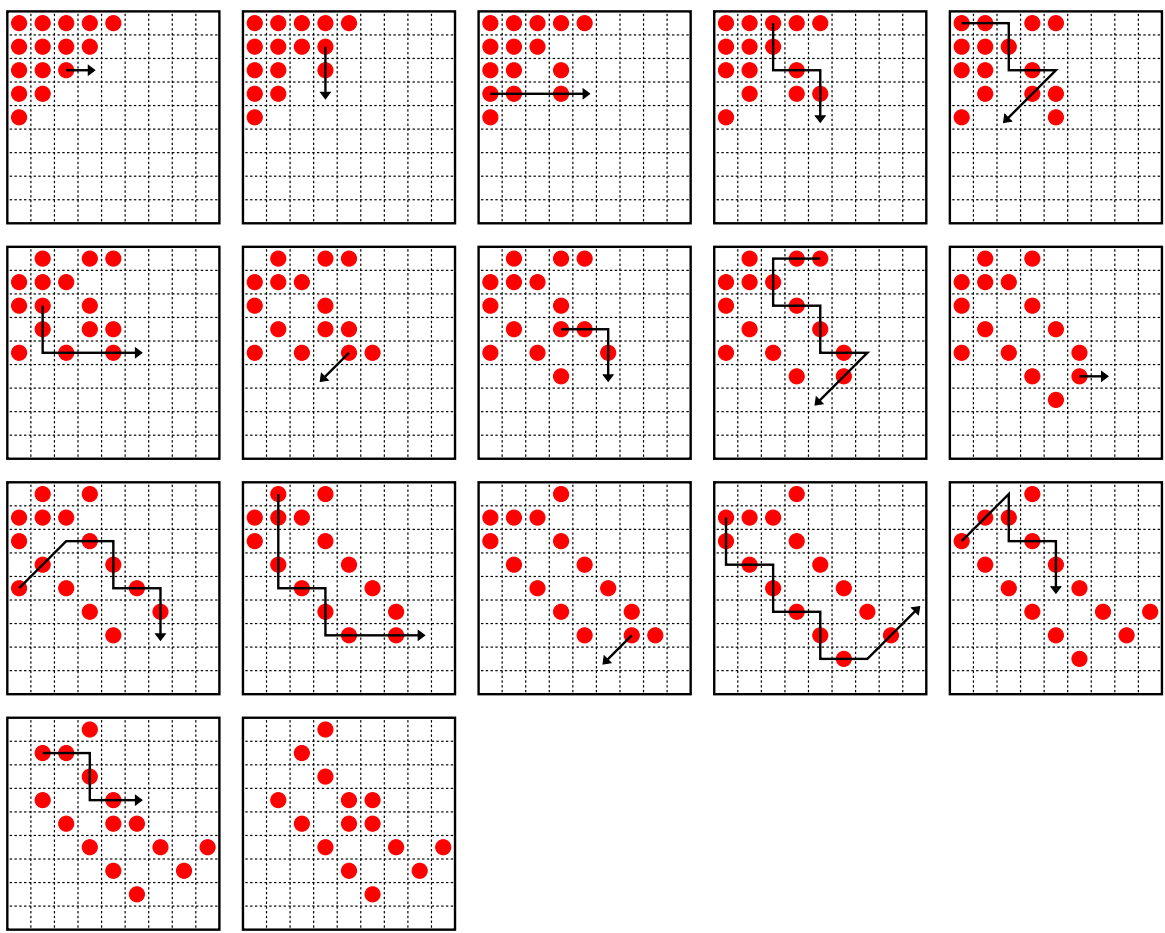}
\caption{The first half of a 31~move solution starting with 15 men.
The only other palindromic 31~move solution is obtained by
extending the 15th move to end at c7.}
\label{fig9}
\end{figure}

Now consider the Chinese Checker transfer problem for the 15-man army
(Problem \#4 in Table~\ref{table3} under 6-move rules).
After $N=15$ moves, for a palindromic solution of length $2N=30$ to exist,
we must have $\mbox{sym}(A_{15})=15$.
The symmetry constraint (\ref{eq:symconst}) is therefore
\begin{equation}
\mbox{sym}(A_i)\ge 15-2(15-i) = 2i-15 \mbox{ for } 8\le i\le 15.
\label{eq:symconst1}
\end{equation}
The symmetry constraint (\ref{eq:symconst1})
can be added to the MCC algorithm,
and it speeds up the search after level 7,
and we can run the search to level 15.
The resulting bidirectional search finds that
there does not exist a \textit{palindromic} 30~move solution.
It remains a possibility that a non-palindromic 30~move
solution exists.
The same search technique finds two palindromic solutions of
length 31 (Figure~\ref{fig9}).

In general, the symmetry constraint (\ref{eq:symconst})
can be added to our search algorithm,
and it can dramatically increase the search speed.
Once again, though, this search technique
contains a ``chicken-in-egg paradox",
because in order to run the search,
we must already know the solution's length.
For a general problem,
our search strategy is first to find a solution of any length,
then try to find a shorter solution, and finally to prove that no shorter
solution exists.

To find a solution of any length, we can truncate the level sets at
each step, as explained in the previous Section~3.2.
This technique can be viewed as a truncation based on the centroid
``score" $c(A)$, but it does not place any extra value on palindromic solutions.
If the truncated $n$th level set is $R_n$, then to find a solution of
length $2N$ the set $R_N$ must satisfy the conditions:
\begin{eqnarray}
\max\left\{c(A), \forall A\in R_N\right\} & \ge & 0, \\
\min\left\{c(A), \forall A\in R_N\right\} & \le & 0.
\end{eqnarray}
Since the truncation technique keeps board positions with largest
centroid, the first condition will eventually be satisfied as $N$ increases.
However, for large problems, the second condition may not be satisfied (because of the truncation),
and then no solution will ever be found, no matter how large $N$ is.
For problems starting from 15-man armies, we instead truncate the level sets
based on the modified score
$$c(A)+\beta \mbox{sym}(A),$$
where $\beta\ge 0$ is an arbitrary weight factor for symmetry.
This modified score weights palindromic solutions more heavily.
We have had good results using $\beta=2$.

The palindromic solutions in Figures~\ref{fig6}--\ref{fig9}
are symmetric about the line $x=y$.
Note that the Checkers board of Figure~\ref{fig2}
is not symmetric about the line $x=y$;
a palindromic solution on this board must instead
be symmetric with respect to $180^\circ$ rotation.
Martin Gardner gives a 20~move palindromic solution on this
board \cite[p. 217]{GardnerNMD}, and our algorithm can find all
palindromic solutions,
using a version of $\mbox{sym}(A)$ corresponding to a
rotation rather than a reflection.

\vskip 30pt
\section*{\normalsize 3.4 Fast armies and balanced armies}
Given an army $A$, we define the center of mass
$\mathbf{m}(A)\in \mathbb{R}^2$ as the average coordinate over the army,
\begin{equation}
\mathbf{m}(A) = \frac{1}{|A|} \sum_{\mathbf{a}\in A} (a_x, a_y).
\label{eq:defcom}
\end{equation}
The diagonal centroid of Section~3.1
is related to $\mathbf{m}(A)=(m_x,m_y)$ by $c(A)=|A|(m_x-m_y)$.
If an army $A$ moves to $B$ in $n$ moves,
we define its average \textbf{speed}
as the distance traveled by the center of mass per move,
\begin{equation}
\sigma = \frac{d(\mathbf{m}(A), \mathbf{m}(B))}{n}=\frac{||\mathbf{m}(A)-\mathbf{m}(B)||}{n},
\label{eq:defspeed}
\end{equation}
using the norm (\ref{eq:normdef}) specific to 4, 6 or 8-move rules.
Note that the same army may have a different speed
measured under different rules.
We call the army $B$ a \textbf{translate} of $A$ if
$B$ can be obtained by shifting each man in $A$ by the
same (2D) vector $\mathbf{k}$,
$B=\{\mathbf{a}+\mathbf{k}, \forall \mathbf{a} \in A\}$.

\begin{figure}[htb]
\centering
\epsfig{file=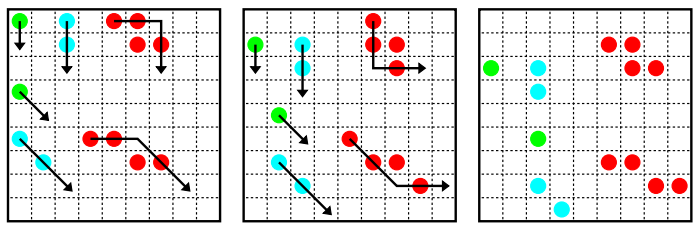}
\caption{Speed 1 armies under 4 and 8-move rules:
``atom" (1 man, green), ``frog" (2 men, blue), ``serpent" (4 men, red).}
\label{fig10}
\end{figure}

Figure~\ref{fig10} shows fast armies of size 1, 2, or 4.
Under 4-move rules,
the three armies in the top row each have speed 1.
Auslander et al. \cite{OptLeaping93} prove that if
$A$ moves to $B$ under 4-move rules,
and $B$ is a translate of $A$,
then the speed $\sigma$ cannot exceed $1$,
and the only speed $1$ armies are the three in the top row
of Figure~\ref{fig10}.
The bottom row of Figure~\ref{fig10} shows armies under 8-move rules,
all having speed $1$.
Under 8-move rules, all six armies in Figure~\ref{fig10}
have speed $1$, except for the top serpent which
has speed $1/2$.

We say that an army is \textbf{balanced} if the
distribution of men over types is as uniform as possible.
The smallest example of a perfectly balanced army are the serpents
in Figure~\ref{fig10},
which contain one man of each type.
The balance of an army can only be changed by a step move.
One reason why a balanced army may be fast is that the
number of possible jump moves is large.
For example, among all 4-man armies under 4-move rules,
the largest number of possible jump moves is 8,
and this can only be achieved by taking one man of each type,
as in the serpent, or 4 men in a square configuration.

In the previous sections we have seen fast armies of larger size.
The speed in each case can be computed and is given in Appendix~A.
One interesting aspect of the armies in
Figures~\ref{fig6}, \ref{fig7} and \ref{fig9}
is that they are usually balanced.
In the Chinese Checker problem of Figure~\ref{fig6},
the army is initially unbalanced,
because there are $(3,3,3,1)$ men of types 0--3.
The first move (d1-d2) converts a man of type 1 to type 3,
resulting in the balanced army $(3,2,3,2)$.

Finally, we note that the most common opening among
experienced players of competitive Chinese Checkers is
the balance restoring step d1-d2 \cite{CCGuide}.
Is this a coincidence?
In the future, it may be productive to look
into the concept of army balance in
competitive Chinese Checkers games.
Halma begins with the unbalanced army $(6,5,5,3)$;
restoring balance would suggest a diagonal step ending at
the type $3$ cells d4, f2, b6 or f6 as one of the opening moves.

\vskip 30pt
\section*{\normalsize 4. Summary}

Chinese Checkers and Halma seem to have a reputation as slow-moving games.
Alternate (but more complex) jumping rules have even been devised
in an attempt to speed them up \cite{NRCC}.
However, increasing complexity does not necessarily improve a game.
We have shown than an entire game of Chinese Checkers can take as few
as 30~moves.
While such short solutions rely on the two players cooperating,
a competitive game can also move along more quickly than many people think.

We have demonstrated, using computational search,
that 30~moves is the shortest possible Chinese Checkers game.
We have also studied a solitaire version of the game where
the goal is to transport a single army across the board as
quickly as possible.
These solutions often have palindromic symmetry,
and pass through a board position which is mirror symmetric
across the line $x=y$.
We have obtained a bound on how quickly the centroid can increase,
and applied it numerically to show that no standard Chinese
Checkers army can cross the board in under 27~moves.

We have also studied the army transfer problem in a more general context,
considering various army configurations under several different jumping rules,
with results summarized in Table~\ref{table3}.
Can we always find a shortest solution that is palindromic?
The answer is probably no, but we have not seen a counterexample.

The problem of the fastest way to transport a 19-man Halma army across
a $16\times 16$ board is difficult to answer computationally.
Nonetheless, this is an interesting problem to work on by hand,
and we have found a solution in 47~moves, given in Appendix A.
Can the reader find a shorter solution?
The centroid arguments of Section~3.1 give a lower bound of only 28~moves,
so there would appear to be ample room for improvement.

Finally, we have introduced the concept of a balanced army,
which means to keep (as much as possible)
the same number of men of each type.
The fastest armies usually seem to be balanced,
and this is an interesting area for further study.

\noindent
\textbf{Acknowledgments}
We thank the anonymous referee for several suggestions which
improved the paper.

\section*{\normalsize Appendix A. Solutions}

\subsection*{\normalsize A.1 Shortest games}

\noindent
Chinese Checkers in 30 moves (Figure~\ref{fig5}), 10 man armies, 6-move rules (by David Fabian):
\begin{small}
\textcolor{blue}{c2-d2}, \textcolor{red}{h8-h6},
\textcolor{blue}{d1-d3}, \textcolor{red}{i6-g6},
\textcolor{blue}{a3-c3-e3}, \textcolor{red}{g9-g7-g5},
\textcolor{blue}{a2-c2-e2-e4}, \textcolor{red}{h7-h5-f7},
\textcolor{blue}{e4-f4}, \textcolor{red}{i9-g9-g7-e7},
\textcolor{blue}{d3-c4}, \textcolor{red}{g6-g4-e4-e2-c2-a2},
\textcolor{blue}{a4-c2-e2-e4-g4-g6-e8-e6}, \textcolor{red}{i8-i6-g6-g4-e4-e2-c2-a4},
\textcolor{blue}{a1-a3-c3-c5}, \textcolor{red}{f9-h7-h5-f5-f3-d3-d1},
\textcolor{blue}{c5-d5}, \textcolor{red}{e7-e5-c5-c3-a3-a1},
\textcolor{blue}{b2-b4-d4-d6-f6-f8-h8}, \textcolor{red}{i7-g9-g7-e7-e5-c5-c3-a3},
\textcolor{blue}{c1-e1-c3-c5-e5-e7-g7}, \textcolor{red}{h6-f8-f6-d6-d4-b4-b2},
\textcolor{blue}{b3-b4}, \textcolor{red}{g8-g6-g4-e4-e2-c2},
\textcolor{blue}{b1-b3-b5-d3-f3-f5}, \textcolor{red}{g5-e5-c5-c3-c1},
\textcolor{blue}{f5-f6}, \textcolor{red}{f7-f5-f3-d3-b5-b3-b1},
\textcolor{blue}{e6-g6-g8-i8}, \textcolor{red}{h9-h7-f7-f5-f3-d3-b5-b3}
(\textcolor{red}{red} wins).
\end{small}

\noindent
Chinese Checkers in 36 moves, 15 man armies, 6-move rules:
\begin{small}
\textcolor{blue}{e1-e2}, \textcolor{red}{g8-g6},
\textcolor{blue}{c1-e1-e3}, \textcolor{red}{h6-f6},
\textcolor{blue}{e3-e4}, \textcolor{red}{f9-f7-f5},
\textcolor{blue}{a1-c1-e1-e3-e5-g5-e7}, \textcolor{red}{g7-g5-e5-e3-e1-c1-a1},
\textcolor{blue}{a5-b5}, \textcolor{red}{i7-g7-g5-e5-e3-e1-c1},
\textcolor{blue}{c3-a5-c5}, \textcolor{red}{g9-i7-g7-g5-e5-e3-e1-c3-a5},
\textcolor{blue}{a3-c3-e1-e3-e5-g5-g7-i7-g9}, \textcolor{red}{i9-i7-g7-g5-e5-e3-e1-c3-a3},
\textcolor{blue}{a4-c4-c6}, \textcolor{red}{e9-g7-g5-e5-e3-e1-c3},
\textcolor{blue}{c6-d6}, \textcolor{red}{i5-i7-g7-g5-e5-e3-e1},
\textcolor{blue}{a2-a4-c4-c6-e6-g4}, \textcolor{red}{i8-g8-e8-e6-c6-c4-a4-a2},
\textcolor{blue}{d2-f2}, \textcolor{red}{i6-g8-e8-e6-c6-c4-a4},
\textcolor{blue}{b3-d3-f1-f3-d5-d7-f7-h5}, \textcolor{red}{f5-f7-d7-d5-f3-f1-d3-b3},
\textcolor{blue}{d1-f1-d3}, \textcolor{red}{h8-h6-h4-f4-d4-d2},
\textcolor{blue}{b1-d1-f1-f3-d5-d7-f7-f5}, \textcolor{red}{h9-f9-f7-d7-d5-f3-f1-d1-b1},
\textcolor{blue}{b4-b6-d4-f4-h4-h6-h8}, \textcolor{red}{h7-h9-f9-f7-d7-d5-f3-f1-d1},
\textcolor{blue}{c2-c4-c6-e6-e8}, \textcolor{red}{g6-e6-c6-c4-c2},
\textcolor{blue}{b2-b4-b6-d4-f4-h4-h6}, \textcolor{red}{f6-f4-d4-b6-b4-b2},
\textcolor{blue}{d3-f1-f3-d5-d7-f7-f9-h9}, \textcolor{red}{f8-d8-f6-f4-d4-b6-b4}
(\textcolor{red}{red} wins).
\end{small}

\subsection*{\normalsize A.2 Short solutions to army transfer problems}

Problem numbers refer to those in Table~\ref{table3}.
All solutions given are palindromic, so only
half the solution is given, followed by ``(reflect)".
This indicates to repeat the moves in reverse order,
reflected about the line $x=y$.
All solutions are the shortest possible, except as noted.

\noindent
Problem \#2 (Figure~\ref{fig6}),
6-move rules in 27 moves (Octave Levenspiel), $C^{max}_{13}=5$, $\sigma=4/9$:
\begin{small}
d1-d2, b1-d1-d3, a3-c3-e3, e3-f3, c1-c3-e3-g3, a1-a3-c3-e3,
c2-e2-e4-g2-g4, g4-g5, a4-c2-e2-e4-g2-g4-g6, a2-c2-e2-e4-g2-g4,
d2-d4-f2-f4-h4-f6-h6, h6-h7, b3-c2,
b2-d2-d4-f2-f4- (reflect).
\end{small}
An alternate solution with only 8 step moves:
\begin{small}
c2-d2, a4-c2-e2, c1-e1-e3, e3-e4, a1-c1-e1-e3-e5,
a3-c1-e1-e3, d1-d3-f3-d5-f5, f5-g5, a2-c2, b3-d1-d3-f3-d5-f5-h5,
b1-b3-d1-d3-f3-d5-f5, d2-f2-d4-f4-f6-h4-h6, h6-h7, b2-d2-f2-d4-f4-
(reflect).
\end{small}

\noindent
Problem \#2, 8-move rules in 20 moves, $C^{max}_{10}=12$, $\sigma=3/10$:
\begin{small}
a2-c4, c4-d4, a4-a2-c4-e4, b1-d3-f5, f5-f6,
a3-c3-e5-g7, d1-b1-d3-f5-f7-h7, g7-h8, a1-c3-e5-g7-i7, b3-b1-d3-f5,
(reflect).
\end{small}

\noindent
Problem \#2, 4-move rules in 30 moves, $C^{max}_{14}=5$, $C^{max}_{15}=11$, $\sigma=2/5$:
\begin{small}
a4-b4, a2-a4-c4, b4-d4, d4-d5, c1-c3-c5-e5, e5-e6, b2-b4-d4-d6-f6, f6-f7,
a3-c3-c5-e5-e7-g7, g7-h7, d1-c1\footnote{Note this backward step move.}, c1-c3-c5-e5-e7-g7-i7, i7-i8, b3-c3,
c3-c5-e5-e7-g7-i7-i9,
(reflect).
\end{small}

\noindent
Problem \#6 (Figure~\ref{fig7}), 8-move rules in 16 moves,
$C^{max}_{8}=10$, $\sigma=3/8$:
\begin{small}
b2-d4, c3-e5, e5-f6, a3-c3-e5-g7, c1-c3-e5,
a1-a3-c3, a2-c4-e4-e6-g6-g8, b1-d3-d5-f5-f7-h7,
(reflect).
\end{small}

\noindent
Problem \#3 (Figure~\ref{fig8}), 6-move rules, jumps only in 35 moves,
$C^{max}_{17}=8$, $\sigma=12/35$:
\begin{small}
b2-d2, d1-d3, c1-c3-e3, a1-c1, c2-e2-c4,
a4-c2-e2-e4, d2-d4-f4, a2-a4-c2, d3-f3-f5, f4-d4-b4-b2-d2,
b1-d1-d3-f3-d5, c1-c3-c5-e5-g5, a3-c3-c5-e5, e4-e6-g4-g6, b3-d1-d3,
c2-e2-e4-e6, f5-h5-f7, d2-d4-d6-
(reflect).
\end{small}

\noindent
Problem \#4 (Figure~\ref{fig9}), 6-move rules in 31 moves,
possibly not shortest, $\sigma=32/93$:
\begin{small}
c3-d3, d2-d4, a4-c4-e4, c1-c3-e3-e5, a1-c1-c3-e3-c5, b3-b5-d5-f5, e5-d6,
d4-f4-f6, e1-c1-c3-e3-e5-g5-e7, f6-g6, a5-c3-e3-e5-g5-g7, b1-b3-b5-d5-d7-f7-h7,
g7-f8, a2-a4-c4-c6-e6-e8-g8-i6, a3-c1-c3-e3-e5,
b2-d2-d4-f4-
(reflect).
\end{small}

\noindent
A 19-man Halma army crosses the board in 47 moves,
probably not shortest, $\sigma=225/893$:
\begin{small}
d2-d4, c3-e5, e5-f6, a1-c3-e5-g7, g7-h8,
c1-c3-e5-g7-i9, i9-j10, e1-c1-c3-e5-g7-i9, a3-c3-e5-g7, a5-c3-e5,
b3-d5-f5-f7-h7-h9-j9-j11, j11-k11, b5-b3-d5-f5-f7-h7-h9-j9-j11-l11, l11-l12, d1-b3-d5-f5-f7-h7-h9-j9-j11-l11-l13,
l13-m13, c4-e4-e6-g6-g8-i8-i10-k10-k12-m12-m14, m14-n14, e2-c4-e4-e6-g6-g8-i8-i10-k10-k12-m12-m14-o14, d3-d5-f5-f7-h7-h9-j9-j11-l11-l13-n13-n15,
a4-c4-e4-e6-g6-g8-i8-i10-k10-k12-\\m12-m14-o16, b2-b3, b4-c3, b1-d3-d5-f5-f7-h7-h9- (reflect).
\end{small}


\end{document}